\documentclass[noinfoline]{imsart}
\usepackage{amsthm,amsmath, amssymb}
\usepackage[square,comma, numbers]{natbib}
\usepackage{graphicx}
\setattribute{journal}{name}{}
\setattribute{journal}{issn}{}
\newcommand{\dist}[1]{\vspace{8pt}\noindent{\bf#1}}
\usepackage{hyperref}

\hypersetup{ 
	pdfauthor={Gavin E. Crooks}, 
	pdftitle={The Amoroso Distribution} 
}

\begin{document}

\title{The Amoroso Distribution}

\author{Gavin E. Crooks}
\address{Physical Biosciences Division,\\ Lawrence Berkeley National Laboratory,\\ Berkeley, CA 94720\\ gecrooks@lbl.gov}

\begin{abstract}
Herein, we review the properties of the Amoroso distribution, the natural unification of the gamma and extreme value  distribution families. Over 50 distinct, named distributions (and twice as many synonyms) occur as special cases or limiting forms. Consequently,  this single simple functional form  encapsulates and systematizes an extensive menagerie of interesting and common probability distributions. 
\end{abstract}

\maketitle


\tableofcontents

\section{The Amoroso  distribution family}
The {\bf Amoroso}  (generalized gamma, Stacy-Mihram) distribution~\cite{Amoroso1925,Johnson1994} is a four parameter,  continuous, univariate, unimodal probability density, with semi-infinite range. The functional form in the most straightforward parameterization is
\begin{align}
\label{Amoroso}  
 \text{Amoroso}(x| a, \theta, \alpha, \beta) 
&=
\frac{1}{\Gamma(\alpha)} 
\left|\frac{\beta}{\theta}\right|
\left(\frac{x-a}{\theta}\right)^{\alpha \beta -1}
\exp \left\{
-  \left(\frac{x-a}{\theta}\right)^{\beta}
\right\}
\\ \notag
& \text{for } x,\ a,\ \theta,\ \alpha,\ \beta\  \text{in } \mathbb{R}, 
\ \alpha>0, \ 
\\ \notag
& \text{support } x \geq a \ \text{if}\ \theta > 0,  \ x\leq a  \ \text{if}\  \theta < 0 .
\end{align}

The Amoroso distribution was originally developed to model lifetimes~\cite{Amoroso1925}. It occurs as the Weibullization of the standard gamma distribution (\ref{Gamma}) and, with integer $\alpha$, in extreme value statistics (\ref{GenFisherTippett}). The Amoroso distribution is itself a limiting form of various more general distributions, most notable the generalized beta and generalized beta prime distributions~\cite{McDonald1984}.

A useful and important property of the Amoroso distribution is that many common and interesting probability distributions occur as  special cases or limits  (See Table~\ref{AmorosoTable}). (Informally, an ``interesting distribution'' is one that has acquired a name, which generally indicates that the distribution is the solution to one or more interesting problems.) This provides a convenient method for systemizing  a significant fraction of the probability distributions that are encountered in practice, provides a consistent parameterization for those distributions, and obviates the need to enumerate the properties (mean, mode, variance, entropy and so on) of each and every specialization.

\paragraph*{Notation:} The four real parameters of the Amoroso distribution consist of a location parameter~$a$, 
a scale parameter~$\theta$,  and two shape parameters,~$\alpha$ and~$\beta$. Whenever these symbols appears in special cases or limiting forms, they refer directly to the parameters of the Amoroso distribution.
The shape parameter $\alpha$ is positive, and in many specializations an integer, $\alpha=n$, or half-integer, $\alpha=\tfrac{k}{2}$. The negation of a standard parameter is indicated by a bar, e.g. $\bar{\beta} = -\beta$. The chi, chi-squared and related distributions are traditionally parameterized with the scale parameter $\sigma$, where $\theta= (2\sigma^2)^{1/{\beta}}$, and $\sigma$ is the standard deviation of a related normal distribution.  Additional alternative parameters are introduced as necessary. 
  
We write  $\text{Amoroso}(x| a, \theta, \alpha, \beta)$ for a density function,  $\text{Amoroso}(a, \theta, \alpha, \beta)$ for the corresponding random variables,  and  $X\sim\text{Amoroso}(a, \theta, \alpha, \beta)$ to indicate that two random variables have the same probability distribution~\cite{Gelman2004}.

\begin{table}[tph!]
\label{AmorosoTable}
\caption{The Amoroso family of distributions.}
\begin{align*}
 \text{Amoroso}(x| a, \theta, \alpha, \beta) 
= &
\frac{1}{\Gamma(\alpha)} 
\left|\frac{\beta}{\theta}\right|
\left(\frac{x-a}{\theta}\right)^{\alpha \beta -1}
\exp \left\{
-  \left(\frac{x-a}{\theta}\right)^{\beta}
\right\}
\\ \notag
& \text{for } x,\ a,\ \theta,\ \alpha,\ \beta\  \text{in } \mathbb{R}, 
\ \alpha>0, \ 
\\ \notag & k,\ n\ \text{positive integers}
\\ \notag
& \text{support } x \geq a \ \text{if}\ \theta > 0,  \ x\leq a  \ \text{if}\  \theta < 0 .
\end{align*}

\begin{tabular}{llccccl}
(\ref{Amoroso}) &Amoroso & $a$ & $\theta$ & $\alpha$ & $\beta$
\\ \hline
(\ref{Stacy}) & Stacy & $0$ & . & . & . \\
(\ref{GenFisherTippett}) & gen. Fisher-Tippett  & . & . & $n$ & .  \\
(\ref{FisherTippett}) &  Fisher-Tippett & . & . & 1 & .  \\
(\ref{Frechet}) &Fr\'{e}chet  & . & . & 1 &  $<\!\!0$  \\
(\ref{GenFrechet}) &  generalized Fr\'{e}chet & . & . & $n$ & $<\!\!0$ \\
(\ref{ScaledInvChi}) &scaled inverse chi& 0 & . & $\tfrac{1}{2}k$  & -2  \\
(\ref{InvChi}) & inverse chi  & 0 & $\frac{1}{\sqrt{2}}$ & $\tfrac{1}{2}k$ & -2 \\
(\ref{InvRayleigh}) &  inverse Rayleigh  & $0$ & . & $1$ & -2 \\
(\ref{PearsonV}) & Pearson type V & . & . & . & -1 \\
(\ref{InvGamma}) & inverse gamma & 0 & . & . & -1 \\
(\ref{ScaledInvChiSqr}) & scaled inverse chi-square & 0 & . & $\tfrac{1}{2}k$ & -1 \\
(\ref{InvChiSqr}) & inverse chi-square & 0 & $\frac{1}{2}$ & $\tfrac{1}{2}k$ & -1 \\
(\ref{Levy}) & L\'{e}vy &  . & . & $\frac{1}{2}$ & -1 \\
(\ref{InvExp}) &  inverse exponential & 0  & . & 1 & -1 \\
(\ref{PearsonIII}) &Pearson type III  & . & . & . & 1 \\
(\ref{Gamma}) &gamma & $0$ & . & . & $1$ \\
(\ref{Gamma}) & Erlang & $0$ & $>\!\!0$ & $n$ & $1$ \\
(\ref{StdGamma}) &standard gamma & 0 & 1 & . & 1  \\ 
(\ref{ScaledChiSqr}) & scaled chi-square & 0 & . & $\tfrac{1}{2}k$ & 1 \\
(\ref{ChiSqr}) & chi-square & 0 & 2 & $\tfrac{1}{2}k$ & 1 \\
(\ref{ShiftExp}) & shifted exponential & . & . & 1 & 1  \\ 
(\ref{Exp}) & exponential & $0$ & . & $1$ & $1$ \\
(\ref{Exp}) &standard exponential & 0 & 1 & 1 & 1 \\ 
(\ref{Gamma}) & Wien & 0 & . & 4& 1 \\
(\ref{Nakagami}) & Nakagami & . & . & . & $2$ \\
(\ref{ScaledChi}) &scaled chi& 0 & . & $\tfrac{1}{2}k$  & 2  \\
(\ref{Chi}) & chi & 0 & $\sqrt{2}$ & $\tfrac{1}{2}k$ & 2 \\
(\ref{HalfNormal}) & half-normal & 0 & . & $\tfrac{1}{2}$ & 2 & \\  
(\ref{Rayleigh}) & Rayleigh & 0 & . & 1 & 2  \\
(\ref{Maxwell}) & Maxwell& 0 & . & $\frac{3}{2}$  & 2  \\
(\ref{WilsonHilferty}) &Wilson-Hilferty& 0 & . & .  & 3  \\
(\ref{GenWeibull}) & generalized Weibull  & . & . & $n$ & $>\!\!0$  \\
(\ref{Weibull}) & Weibull & . & . & 1 &  $>\!\!0$  \\
(\ref{PseudoWeibull}) & pseudo-Weibull & . & . & $1$+$\tfrac{1}{\beta}$ &  $>\!\!0$  \\
(\ref{StretchedExp}) & stretched exponential & 0 & . & 1 & $>\!\!0$ & \\

\\
 & \underline{Limits} \\
(\ref{LogGamma}) & log-gamma  &. & . & . & . &$\lim_{\beta\rightarrow\infty}$  \\
(\ref{PowerLaw}) & power law &.&.& $\tfrac{1-p}{\beta}$& . &$\lim_{\beta\rightarrow0}$  \\
(\ref{LogNormal}) & log-normal  & . & . &$\tfrac{1}{(\beta\sigma)^2}$  & . &$\lim_{\beta\rightarrow0}$  \\
(\ref{Normal}) & normal  & . &  . & . & 1 &$\lim_{\alpha\rightarrow\infty}$  \\
\\
\\
\\
\\
\\
\end{tabular} 

\end{table}


\subsection{Special cases: Miscellaneous}

\dist{Stacy} (hyper gamma, generalized Weibull, Nukiyama-Tanasawa, generalized gamma, generalized semi-normal, hydrograph, Leonard hydrograph, transformed gamma)  distribution~\cite{Stacy1962,Dadpay2007}:
\begin{align}
\label{Stacy}
\text{Stacy}(x | \theta, \alpha, \beta) 
=& \frac{1}{\Gamma(\alpha)} \left|\frac{\beta}{\theta}\right| \left(\frac{x}{\theta}\right)^{\alpha\beta-1} 
\exp \left\{ -\left(\frac{x}{\theta}\right)^{\beta} \right\}
\\=&  \text{Amoroso}(x| 0, \theta, \alpha, \beta) \notag
\end{align}
If we drop the location parameter from $\text{Amoroso}$, then we obtain the 
Stacy, or generalized gamma distribution, the parent of the gamma family of distributions.
If $\beta$ is negative then the distribution is  {\bf generalized inverse gamma}, the parent of various inverse distributions, including the inverse gamma (\ref{InvGamma}) and inverse chi (\ref{InvChi}). 

The Stacy distribution is obtained as the positive even powers,  modulus, and powers of the modulus of a centered, normal random variable (\ref{Normal}), 
\[
\text{Stacy}\left((2\sigma^2)^{\tfrac{1}{\beta}} ,\tfrac{1}{2}, \beta\right) \sim \Big|\text{Normal}(0,\sigma)\Big|^{\tfrac{2}{\beta}}
\]
and as powers of the sum of squares of $k$ centered, normal random variables. 
\[
\text{Stacy}\left( (2\sigma^2)^{\tfrac{1}{\beta}} ,\tfrac{1}{2}k, \beta\right) \sim  \left( \sum_{i=1}^{k} \Big(\text{Normal}(0,\sigma)\Big)^2\right)^{\tfrac{1}{\beta}}
\]

\dist{Stretched exponential} distribution~\cite{Laherrere1998}:
\begin{align}
\label{StretchedExp}
\text{StretchedExp}(x|\theta,\beta) = & 
\frac{\beta}{|\theta|}
\left(\frac{x}{\theta}\right)^{\beta -1}
\exp\left\{
-  \left(\frac{x}{\theta}\right)^{\beta}
\right\}
\\ \notag & \text {for } \beta>0
\\= & \text{Weibull}(x| 0, \theta, \beta) \notag  
 \\= & \text{Amoroso}(x| 0, \theta,1, \beta) \notag 
\end{align}
Stretched exponentials are an alternative to power laws for modeling fat tailed distributions. For $\beta=1$ we recover the exponential distribution (\ref{Exp}), and $\beta=0$ a power law distribution  (\ref{PowerLaw}). 

\begin{figure}[t]
\begin{center}
\includegraphics[scale=1]{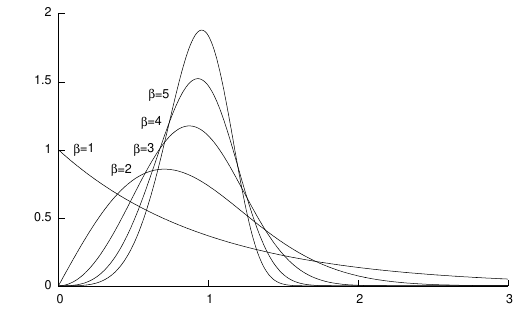}
\end{center}
\caption{$\text{\rm Amoroso}(x|0,1,1,\beta)$, stretched exponential}
\end{figure}

\dist{Pseudo-Weibull} distribution~\cite{Voda1989}:
\begin{align}
\label{PseudoWeibull}
\text{PseudoWeibull}(x|\theta,\beta)  
=& \frac{1}{\Gamma(1+\tfrac{1}{\beta})} \frac{\beta}{|\theta|} \left(\frac{x}{\theta}\right)^{\beta} 
\exp \left\{ -\left(\frac{x}{\theta}\right)^{\beta} \right\}
\\ & \text{for } \beta>0 \notag 
\\=&  \text{Amoroso}(x| 0, \theta, 1+\tfrac{1}{\beta}, \beta) \notag
\end{align}
Proposed as another model of failure times.

\subsection{Special cases: Positive integer $\beta$}

\begin{figure}[t]
\begin{center}
\includegraphics[scale=1]{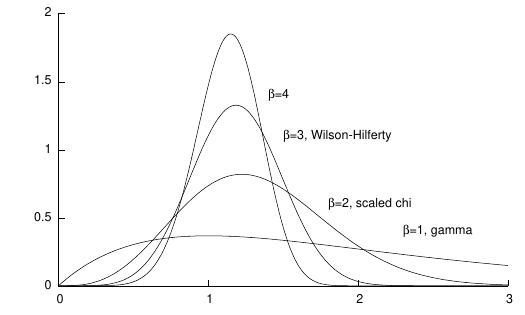}
\end{center}
\caption{$\text{\rm Amoroso}(x|0,1,2,\beta)$}
\end{figure}

\dist{Gamma} ($\Gamma$)  distribution~\cite{Pearson1893, Pearson1895, Johnson1994} : 
\begin{align}
\label{Gamma}
\text{Gamma}(x | \theta, \alpha) 
&=  \frac{1}{\Gamma(\alpha)|\theta|} \left(\frac{x}{\theta}\right)^{\alpha-1} \exp\left\{-\frac{x}{\theta}\right\} \\
&=\text{PearsonIII}(x | 0, \theta, \alpha) \notag \\
&=\text{Stacy}(x| \theta, \alpha,1) \notag 
\\&=  \text{Amoroso}(x| 0, \theta, \alpha, 1) \notag 
\end{align}
The name of this distribution derives from the normalization constant.
The gamma distribution often appear as a solution to problems in statistical physics. For example, the energy density of a classical ideal gas; or the {\bf Wien} (Vienna) distribution $\text{Wien}(x|T)=\text{Gamma}(x|T,4)$, an approximation to the relative intensity of black body radiation as a function of the frequency. The {\bf Erlang} (m-Erlang) distribution~\cite{Erlang1909} is a gamma distribution with integer $\alpha$, which models the waiting time to observe $\alpha$ events from a Poisson process with rate $1/\theta$ ($\theta>0$).

Gamma distributions obey an addition property:
\begin{align*}
\text{Gamma}(\theta, \alpha_1) +  \text{Gamma}(\theta,\alpha_2) \sim \text{Gamma}(\theta,\alpha_1+\alpha_2)
\end{align*}
The  sum of two independent, gamma distributed random variables (with common $\theta$'s, but possibly different $\alpha$'s) is again a gamma random variable~\cite{Johnson1994}.

\begin{figure}[t]
\begin{center}
\includegraphics[scale=1]{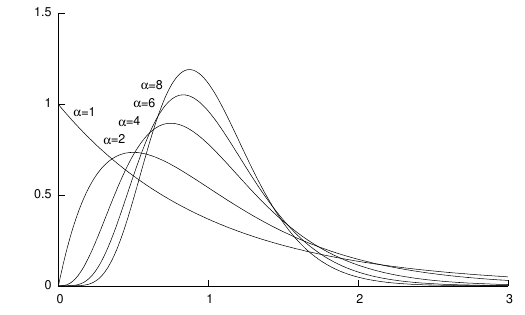}
\end{center}
\caption{$\text{\rm Gamma}(x|\tfrac{1}{\alpha},\alpha)$ (unit variance) }
\end{figure}

\dist{Standard gamma} (standard Amoroso) distribution~\cite{Johnson1994}: 
\begin{equation}
\text{StdGamma}(x|\alpha) = \frac{1}{\Gamma(\alpha)} x^{\alpha-1} e^{-x}
\label{StdGamma}
\end{equation}
The Amoroso distribution can be obtained from the standard gamma distribution by the Weibull change of variables, $x \mapsto \left(\tfrac{x-a}{\theta}\right)^\beta$.
\[
\text{Amoroso}(a ,\theta,\alpha,\beta) \sim
a+\theta \Big[{\text{StdGamma}}(\alpha)\Big]^{1/\beta} 
\]

\dist{Pearson type III} \noindent distribution~\cite{Pearson1895, Johnson1994}:
\begin{align}
\label{PearsonIII}
\text{PearsonIII}(x | a , \theta, \alpha) 
=& \frac{1}{\Gamma(\alpha)|\theta|} \left(\frac{x-a }{\theta}\right)^{\alpha-1} 
\exp\left\{-\left(\frac{x-a }{\theta}\right)\right\}
\\
=&  \text{Amoroso}(x| a  , \theta, \alpha, 1) \notag
\end{align}
The gamma distribution with a location parameter.

\dist{Exponential} (Pearson type X, waiting time, negative exponential) distribution~\cite{Johnson1994}:
\begin{align}
\label{Exp}
\text{Exp}(x | \theta) 
&= \frac{1}{|\theta|} \exp\left\{-\frac{x}{\theta}\right\} \\
&=  \text{Gamma}(x| \theta, 1) \notag \\
&=  \text{Amoroso}(x| 0, \theta, 1, 1)  \notag
\end{align}
An important property of the exponential distribution is that it is memoryless: the conditional probability given that $x>c$, where $c$ is a positive content, is again an exponential distribution with the same scale parameter. The only other distribution with this property is the geometric distribution~\cite{Evans2000}, the discrete analog of the exponential distribution. With $\theta=1$ we obtain a {\bf standard exponential} distribution. See also shifted exponential (\ref{ShiftExp}), stretched exponential (\ref{StretchedExp}) and inverse exponential (\ref{InvExp}).

\dist{Shifted exponential} distribution~\cite{Johnson1994}:
\begin{align}
\label{ShiftExp}
\text{ShiftExp}(x |a , \theta) 
&=  \frac{1}{|\theta|}  \exp\left\{-\left(\frac{x-a }{\theta}\right)\right\}\\
&= \text{PearsonIII}(x | a , \theta, 1) \notag  \\
&=  \text{Amoroso}(x| a ,\theta,1,1)  \notag 
\end{align}
The exponential distribution with a location parameter.

\dist{Nakagami} (generalized normal, Nakagami-m) distribution~\cite{Nakagami1960}:
\begin{align}
\label{Nakagami}
 \text{Nakagami}(x | a , \theta, m) 
& =
 \frac{2}{\Gamma(\tfrac{m}{2}) |\theta| }
\left(\frac{x-a }{\theta}\right)^{m -1}
\exp \left\{
-  \left(\frac{x-a }{\theta}\right)^{2}
\right\}
\\ \notag
& = \text{Amoroso}(x|a,\theta, \tfrac{m}{2} ,2)
\notag
\end{align}
Used to model attenuation of radio signals that reach a receiver by multiple paths~\cite{Nakagami1960}.

\dist{Half-normal} (semi-normal, positive definite normal, one-sided normal) distribution~\cite{Johnson1994}:
\begin{align}
\label{HalfNormal}
\text{HalfNormal}(x | \sigma ) 
&= \frac{2}{\sqrt{2\pi \sigma^2}} 
\exp\left\{-\left( \frac{x^2}{2\sigma^2}\right) \right\} \\
&=\text{ScaledChi}(x |  \sigma, 1) \notag \\
&=  \text{Stacy}(x| \sqrt{2\sigma^2} ,\tfrac{1}{2},2) \notag \\
&=  \text{Amoroso}(x| 0, \sqrt{2\sigma^2} , \tfrac{1}{2}, 2) \notag 
\end{align}
The modulus of a normal distribution with zero mean and variance $\sigma^2$.

\dist{Chi-square} ($\chi^2$)  distribution~\cite{Fisher1924,Johnson1994}:
\begin{align}
\label{ChiSqr}
\text{ChiSqr}(x | k) 
&= \frac{1}{2\Gamma(\tfrac{k}{2})} \left(\frac{x}{2}\right)^{\tfrac{k}{2}-1} 
\exp\left\{-\left(\frac{x}{2}\right)\right\} \\
& \qquad \text{for positive integer } k \notag \\
&=  \text{Gamma}(x| 2,\tfrac{k}{2}) \notag \\
&=  \text{Stacy}(x|2, \tfrac{k}{2},1) \notag \\
&=  \text{Amoroso}(x| 0, 2, \tfrac{k}{2}, 1) \notag 
\end{align}
The distribution of a sum of squares of $k$ independent standard normal random variables.  The chi-square distribution is important for statistical hypothesis testing in the frequentist approach to statistical inference.

\begin{figure}[t]
\begin{center}
\includegraphics[scale=1]{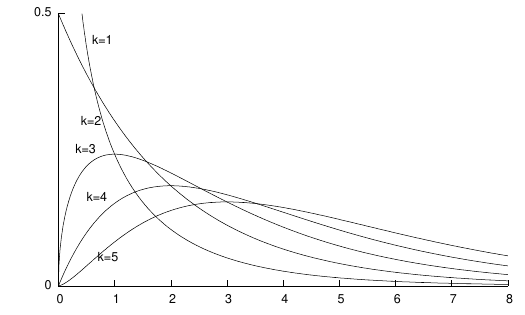}
\end{center}
\caption{$\text{\rm ChiSqr}(x|k)$}
\end{figure}

\dist{Scaled chi-square}  distribution~\cite{Lee2012}:
\begin{align}
\label{ScaledChiSqr}
\text{ScaledChiSqr}(x | \sigma, k) 
&= \frac{1}{2\sigma^2\Gamma(\tfrac{k}{2})} \left(\frac{x}{2\sigma^2}\right)^{\tfrac{k}{2}-1} 
\exp\left\{-\left(\frac{x}{2\sigma^2} \right)\right\} \\
& \qquad \text{for positive integer } k \notag \\
&=  \text{Stacy}(x|2\sigma^2, \tfrac{k}{2},1) \notag \\
&=\text{Gamma}(x|2\sigma^2, \tfrac{k}{2}) \notag \\
&=  \text{Amoroso}(x| 0, 2\sigma^2, \tfrac{k}{2}, 1) \notag 
\end{align}
The distribution of a sum of squares of $k$ independent normal random variables with variance $\sigma^2$.

\dist{Chi} ($\chi$) distribution~\cite{Johnson1994}:
\begin{align}
\label{Chi}
\text{Chi}(x | k) 
&= \frac{ \sqrt{2}}{\Gamma(\tfrac{k}{2})} { \left(\frac{x}{\sqrt{2}}\right)}^{k-1} 
\exp\left\{ -\left( \frac{x^2}{2}    \right)\right\} 
\\
& \qquad \text{for positive integer } k \notag \\
& = \text{ScaledChi}(x|1,k) \notag \\
&=  \text{Stacy}(x|\sqrt{2}, \tfrac{k}{2}, 2)  \notag \\
&=  \text{Amoroso}(x| 0, \sqrt{2} , \tfrac{k}{2}, 2) \notag 
\end{align}
The root-mean-square of $k$ independent standard normal variables, or the square root of a chi-square random variable.
\[
\text{Chi}(k) \sim \sqrt{\text{ChiSqr}(k)}
\]

\dist{Scaled chi} (generalized Rayleigh) distribution~\cite{Miller1964,Johnson1994}:
\begin{align}
\text{ScaledChi}(x | \sigma, k) 
&= \frac{2}{\Gamma(\tfrac{k}{2}) \sqrt{2\sigma^2}} { \left(\frac{x}{\sqrt{2\sigma^2}}\right)}^{k-1} 
\exp\left\{-\left(\frac{x^2}{2\sigma^2}\right)\right\} 
\notag
\\
& \qquad \text{for positive integer } k \notag \\
&=  \text{Stacy}(x|\sqrt{2\sigma^2}, \tfrac{k}{2},2) 
\label{ScaledChi}
\\
&=  \text{Amoroso}(x | 0, \sqrt{2\sigma^2}, \tfrac{k}{2}, 2) 
\notag 
\end{align}
The root-mean-square of $k$ independent and identically distributed normal variables with zero mean and variance~$\sigma^2$.

\begin{figure}[t]
\begin{center}
\includegraphics[scale=1]{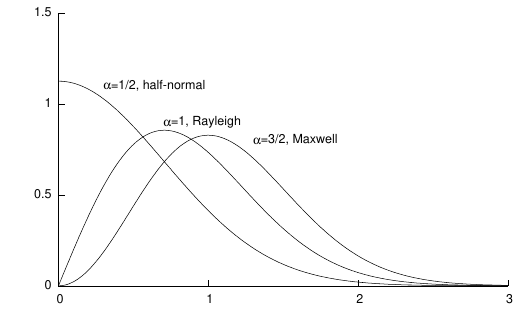}
\end{center}
\caption{$\text{\rm Amoroso}(x|0,1,\alpha,2)$}
\end{figure}

\dist{Rayleigh} distribution~\cite{Strutt1880,Johnson1994}:
\begin{align}
\label{Rayleigh}
\text{Rayleigh}(x | \sigma) 
&= \frac{1}{\sigma^2 }\ x\  \exp\left\{-\left(\frac{x^2}{2 \sigma^2}\right)\right\} 
\\
&=\text{ScaledChi}(x | \sigma, 2) \notag \\
&=  \text{Stacy}(x| \sqrt{2\sigma^2} ,1,2) \notag \\
&=  \text{Amoroso}(x| 0, \sqrt{2\sigma^2} , 1, 2) \notag 
\end{align}
 The root-mean-square of two independent and identically distributed normal variables with zero mean and variance $\sigma^2$. 
 For instance, wind speeds are approximately Rayleigh distributed, since the horizontal components of the velocity are approximately normal, and the vertical component is typically small~\cite{Justus1978}.

\dist{Maxwell} (Maxwell-Boltzmann, Maxwell speed) distribution~\cite{Maxwell1860, Abramowitz1965}:
\begin{align}
\label{Maxwell}
\text{Maxwell}(x | \sigma) 
&= \frac{\sqrt{2}}{\sqrt{\pi} \sigma^3}\ x^2 \exp\left\{-\left(\frac{x^2}{2\sigma^2}\right)\right\}  
 \\
&=\text{ScaledChi}(x | \sigma, 3) \notag \\
&=  \text{Stacy}(x| \sqrt{2\sigma^2} ,\tfrac{3}{2},2) \notag \\
&=  \text{Amoroso}(x| 0, \sqrt{2\sigma^2} , \tfrac{3}{2}, 2) \notag 
\end{align}
The speed distribution of molecules in thermal equilibrium. The root-mean-square of three independent and identically distributed normal variables with zero mean and variance $\sigma^2$.

\dist{Wilson-Hilferty} distribution~\cite{Wilson1931,Johnson1994}:
\begin{align}
\label{WilsonHilferty}
\text{WilsonHilferty}(x | \theta, \alpha) 
&= \frac{3}{\Gamma(\alpha)|\theta|} \left(\frac{x}{\theta}\right)^{3 \alpha-1} \exp\left\{-\left(\frac{x}{\theta}\right)^{3}\right\}
\\ 
&=  \text{Stacy}(x|\theta, \alpha, 3) 
\notag 
\\ &=  \text{Amoroso}(x| 0, \theta, \alpha, 3) \notag
\end{align}
The cube root of a gamma variable follows the Wilson-Hilferty distribution~\cite{Wilson1931}, which has been used to approximate a normal distribution if $\alpha$ is not too small.
A related approximation using quartic roots of gamma variables~\cite{Hawkins1986} leads to   $\text{Amoroso}(x| 0, \theta, \alpha, 4)$.


\subsection{Special cases: Negative integer $\beta$}

\dist{Pearson type V} distribution~\cite{Pearson1901}:
\begin{align}
\label{PearsonV}
\text{PearsonV}(x |a, \theta, \alpha) 
=&
\frac{1}{\Gamma(\alpha)\left|{\theta}\right|} 
\left(\frac{\theta}{x-a }\right)^{\alpha +1}
 \exp\left\{-\left( \frac{\theta}{x-a }
\right)\right\} 
\notag 
\\
=&  \text{Amoroso}(x| a , \theta, \alpha, -1) 
\end{align}
With negative $\beta$ we obtain various ``inverse'' distributions related to distributions with positive $\beta$ by the reciprocal transformation $ (\tfrac{x-a}{\theta} ) \mapsto (\tfrac{\theta}{x-a} )$.
Pearson's type V is the inverse of Pearson's type III distribution.

\dist{Inverse gamma} (Vinci) distribution~\cite{Johnson1994}:
\begin{align}
\label{InvGamma}
\text{InvGamma}(x | \theta, \alpha) 
&= \frac{1}{\Gamma(\alpha) |\theta|} \left(\frac{\theta}{x}\right)^{\alpha+1} 
 \exp\left\{-\left( \frac{\theta}{x}   \right)\right\} 
\\
&=  \text{Stacy}(x| \theta, \alpha,-1)\notag \\
& = \text{PearsonV}(x|0,\theta,\alpha)  \notag \\
&=  \text{Amoroso}(x| 0, \theta, \alpha, -1) \notag 
\end{align}
Occurs as the conjugate prior for an exponential distribution's scale parameter~\cite{Johnson1994}, or the prior for variance of a  normal distribution with known mean~\cite{Gelman2004}.

\dist{Inverse exponential} distribution~\cite{Kleiber2003}:
\begin{align}
\label{InvExp}
\text{InvExp}(x | \theta) 
&= \frac{|\theta|}{x^2 }  \exp\left\{-\left( \frac{\theta}{x}   \right)\right\}   \\
&=  \text{InvGamma}(x| \theta, 1) \notag \\
&=  \text{Stacy}(x|\theta,1,-1) \notag \\
&=  \text{Amoroso}(x| 0, \theta, 1, -1) \notag 
\end{align}
Note that the name ``inverse exponential'' is occasionally used for the ordinary exponential distribution (\ref{Exp}).

\begin{figure}[t]
\begin{center}
\includegraphics[scale=1]{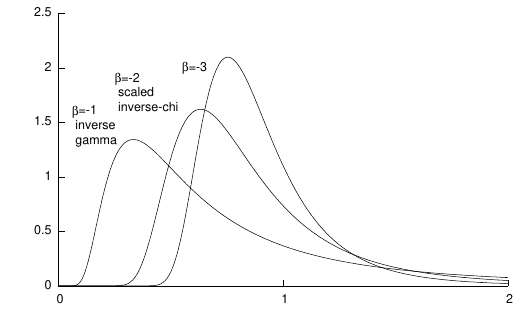}
\end{center}
\caption{$\text{\rm Amoroso}(x|0,1,2,\beta)$, negative $\beta$.}
\end{figure}

\dist{L\'{e}vy} distribution (van der Waals profile)~\cite{Feller1971}: 
\begin{align}
\label{Levy}
\text{L\'{e}vy}(x | a, c) 
&= \sqrt{\frac{c}{2\pi}} \frac{1}{(x-a)^{3/2}}  \exp\left\{-\frac{c}{2(x-a)}\right\}  
\\
&= \text{PearsonV}(x|a,\tfrac{c}{2},\tfrac{1}{2})  \notag \\
&=  \text{Amoroso}(x| a, \tfrac{c}{2}, \tfrac{1}{2}, -1) \notag 
\end{align}
The L\'{e}vy distribution is notable for being stable:  a linear combination of identically distributed  L\'{e}vy distributions is again a  L\'{e}vy distribution. The other stable distributions with analytic forms are the normal distribution (\ref{Normal}), which is also a limit of the Amoroso distribution, and the Cauchy distribution~\cite{Johnson1994}, which is not. L\'{e}vy distributions describe first passage times in one dimensional Brownian diffusion~\cite{Feller1971}.

\dist{Scaled inverse chi-square}  distribution~\cite{Gelman2004}:
\begin{align}
\label{ScaledInvChiSqr}
\text{ScaledInvChiSqr}(x | \sigma, k) 
=& \frac{2 \sigma^2}{\Gamma(\tfrac{k}{2}) } \left(\frac{1}{2 \sigma^2x}\right)^{\frac{k}{2}+1} 
\exp\left\{-\left( \frac{1}{2 \sigma^2x}   \right)\right\}
\\
&\qquad  \text{for positive integer } k \notag \\
&=  \text{InvGamma}(x| \tfrac{1}{2 \sigma^2}, \tfrac{k}{2}) \notag \\
&= \text{PearsonV}(x|0,\tfrac{1}{2 \sigma ^2},\tfrac{k}{2})  \notag  \\
&= \text{Stacy}(x| \tfrac{1}{2 \sigma ^2},\tfrac{k}{2}, -1)  \notag  \\
&=  \text{Amoroso}(x| 0, \tfrac{1}{2 \sigma ^2}, \tfrac{k}{2}, -1) \notag 
\end{align}
A special case of the inverse gamma distribution with half-integer $\alpha$. Used as a prior for variance parameters in normal models~\cite{Gelman2004}.

\dist{Inverse chi-square} distribution~\cite{Gelman2004}: 
\begin{align}
\label{InvChiSqr}
\text{InvChiSqr}(x | k) 
=& \frac{2}{\Gamma(\tfrac{k}{2}) } \left(\frac{1}{2x}\right)^{\frac{k}{2}+1} \exp\left\{-\left( \frac{1}{2x}  \right)\right\}
  \\
&\qquad  \text{for positive integer } k \notag \\
& = \text{ScaledInvChiSqr}(x|1,k)\notag \\
&=  \text{InvGamma}(x| \tfrac{1}{2}, \tfrac{k}{2}) \notag \\
&= \text{PearsonV}(x|0,\tfrac{1}{2},\tfrac{k}{2})  \notag \\
&=  \text{Stacy}(x|\tfrac{1}{2}, \tfrac{k}{2},-1) \notag \\
&=  \text{Amoroso}(x| 0, \tfrac{1}{2}, \tfrac{k}{2}, -1) \notag 
\end{align}
A standard scaled inverse chi-square distribution.

\dist{Scaled inverse chi} distribution~\cite{Lee2012}:
\begin{align}
\label{ScaledInvChi}
\text{ScaledInvChi}(x | \sigma, k) 
&= \frac{2 \sqrt{2 \sigma ^2} }{ \Gamma(\tfrac{k}{2})} { \left(\frac{1}{\sqrt{2 \sigma^2} x}\right)}^{k+1} \exp\left\{-\left(\frac{1}{2 \sigma^2 x^2}  \right)\right\}
\\
&=  \text{Stacy}(x|\tfrac{1}{\sqrt{2 \sigma^2}}, \tfrac{k}{2}, -2)  \notag \\
&=  \text{Amoroso}(x| 0, \tfrac{1}{\sqrt{2 \sigma^2}}, \tfrac{k}{2}, -2) \notag 
\end{align}
Used as a prior for the standard deviation of a normal distribution.

\dist{Inverse chi} distribution~\cite{Lee2012}: 
\begin{align}
\label{InvChi}
\text{InvChi}(x | k) 
&= \frac{2\sqrt{2} }{ \Gamma(\tfrac{k}{2})} { \left(\frac{1}{\sqrt{2} x}\right)}^{k+1} \exp\left\{-\left(\frac{1}{2 x^2}  \right)\right\}
\\
&=  \text{Stacy}(x| \tfrac{1}{\sqrt{2}}, \tfrac{k}{2}, -2)  \notag \\
&=  \text{Amoroso}(x| 0, \tfrac{1}{\sqrt{2}} , \tfrac{k}{2}, -2) \notag 
\end{align}
The standard inverse chi distribution.

\dist{Inverse Rayleigh} distribution~\cite{Evans2000}:
\begin{align}
\label{InvRayleigh}
\text{InvRayleigh}(x | \sigma) 
&= 2 \sqrt{2 \sigma ^2}    \left(\frac{1}{\sqrt{2 \sigma^2} x}\right)^{3} \exp\left\{-\left(\frac{1}{2 \sigma^2 x^2}  \right)\right\}
\\
&=  \text{Stacy}(x|\tfrac{1}{\sqrt{2 \sigma^2}}, 1, -2)  \notag \\
&=  \text{Amoroso}(x| 0, \tfrac{1}{\sqrt{2 \sigma^2}}, 1, -2) \notag 
\end{align}
The inverse Rayleigh distribution has been used to model  failure time~\cite{Voda1972}.


\subsection{Special cases: Extreme order statistics}

\begin{figure}[t]
\begin{center}
\includegraphics[scale=1]{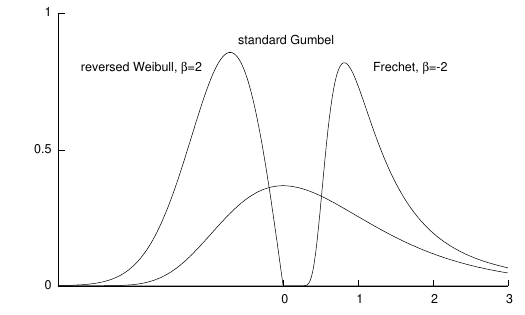}
\end{center}
\caption{Extreme value distributions}
\end{figure}

\dist{Generalized Fisher-Tippett} distribution~\cite{Smirnov1949,Barndorff-Nielsen1963}:
\begin{align}
\label{GenFisherTippett}  
 \text{GenFisherTippett}(x| a, \omega, n, \beta) 
\notag
&=
\frac{n^n}{\Gamma(n)} 
\left|\frac{\beta}{\omega}\right|
\left(\frac{x-a}{\omega}\right)^{n \beta -1}
\exp \left\{
-  n \left(\frac{x-a}{\omega}\right)^{\beta}
\right\}
\\
& \quad \text{for positive integer } n
\\ \notag
& = \text{Amoroso}(x|a,{\omega}/{n^{\frac{1}{\beta} }},n,\beta)
\end{align}
If we take $N$ samples from a probability distribution, then asymptotically for large $N$ and $n\ll N$, the distribution of the $n$th largest (or smallest) sample follows a generalized Fisher-Tippett distribution. The parameter $\beta$ depends on the tail behavior of the sampled distribution. Roughly speaking, if the tail is unbounded and decays exponentially then $\beta$ limits to $\infty$, if the tail scales as a power law then $\beta<0$,  and if the tail is finite $\beta>0$~\cite{Gumbel1958}. In these three limits we obtain the Gumbel (\ref{Gumbel}, \ref{GenGumbel}), Fr\'{e}chet (\ref{Frechet}, \ref{GenFrechet}) and Weibull (\ref{Weibull},\ref{GenWeibull}) families of extreme value distribution (Extreme value distributions types I, II and III) respectively. If $\beta/\omega$ is negative we obtain distributions for the $n$th maxima, if positive then the $n$th minima.

\dist{Fisher-Tippett} (Generalized extreme value, GEV, von Mises-Jenkinson, von Mises extreme value) distribution~\cite{Fisher1928, Mises1936, Gumbel1958,Johnson1995}: 
\begin{align}
\label{FisherTippett}
 \text{FisherTippett}(x| a , \omega, \beta) 
&=
\left|\frac{\beta}{\omega}\right|
\left(\frac{x-a}{\omega}\right)^{ \beta -1}
\exp \left\{
-  \left(\frac{x-a}{\omega}\right)^{\beta}
\right\}
\\ \notag & = \text{GenFisherTippett}(x|a, \omega, 1, \beta)
\\ \notag & = \text{Amoroso}(x|a, \omega, 1, \beta)
\end{align}
The asymptotic distribution of the extreme value from a large sample. The superclass of type I, II and III (Gumbel, Fr\'{e}chet, Weibull) extreme value distributions~\cite{Mises1936}.  This is the distribution for maximum values with $\beta/\omega<0$ and minimum values for $\beta/\omega>0$.

The maximum of two Fisher-Tippett random variables (minimum if $\beta/\omega>0$)  is again a Fisher-Tippett random variable. 
\begin{align*}
\max\Big[ \text{FisherTippett}(a,\omega_1,\beta),  \text{FisherTippett}(a, \omega_2,\beta)  \Big]&\\  \sim 
 \text{FisherTippett}(a, \frac{(\omega_1^{\beta} + \omega_2^{\beta} )^{1/\beta}}{\omega_1 \omega_2},\beta) 
\end{align*}
This follows because taking the maximum of two random variables is equivalent to multiplying their cumulative distribution functions, and the Fisher-Tippett cumulative distribution function is $\exp \left\{
-  \left(\frac{x-a}{\omega}\right)^{\beta}
\right\}$.

\dist{Generalized Weibull} distribution~\cite{Smirnov1949,Barndorff-Nielsen1963}:
\begin{align}
\label{GenWeibull}
\text{GenWeibull}(x | a , \omega, n, \beta) 
&=	\frac{n^n}{\Gamma(n)}  \frac{ \beta}{| \omega |} \left(\frac{x-a }{\omega}\right)^{n \beta-1} \exp\left\{-n \left(\frac{x-a }{\omega}\right)^{ \beta}\right\} 
\\ \notag &\quad \text{for } \beta>0 
\\ \notag & = \text{GenFisherTippett}(x|a, \omega, n, \beta)
\\ \notag
&= \text{Amoroso}(x| a , {\omega}/{n^{\frac{1}{\beta} }}, n, \beta) 
\end{align}
The limiting distribution of the $n$th smallest value of a large number of  identically distributed random variables that are at least~$a$. 
If $\omega$ is negative we obtain the distribution of the $n$th largest value.

\dist{Weibull} (Fisher-Tippett type III, Gumbel type III, Rosin-Rammler, Rosin-Rammler-Weibull, extreme value type III, Weibull-Gnedenko) distribution~\cite{Weibull1951,Johnson1995}: 
\begin{align}
\label{Weibull}
\text{Weibull}(x | a ,\omega, \beta) 
&=	\frac{\beta}{| \omega |} \left(\frac{x-a }{\omega}\right)^{\beta-1} \exp\left\{-\left(\frac{x-a }{\omega}\right)^{\beta}\right\} 
\\ \notag &\quad \text{for } \beta>0 
\\ \notag
& = \text{FisherTippett}(x| a, \omega, \beta) 
\\ \notag
&= \text{Amoroso}(x| a , \omega, 1, \beta) 
\end{align}
This is the limiting distribution of the minimum of a large number of  identically distributed random variables that are at least~$a$.  If $\omega$ is negative we obtain a {\bf reversed Weibull} (extreme value type III) distribution for maxima.
Special cases of the Weibull distribution include the exponential ($\beta=1$) and Rayleigh ($\beta=2$)  distributions.

\dist{Generalized Fr\'{e}chet} distribution~\cite{Smirnov1949,Barndorff-Nielsen1963}:
\begin{align}
\label{GenFrechet}
\text{GenFr\'{e}chet}(x | a , \omega, n, \bar{\beta}) 
 \notag 
&=	\frac{n^n}{\Gamma(n)}  \frac{\bar{\beta}}{| \omega |} \left(\frac{x-a }{\omega}\right)^{-n\bar{\beta}-1} 
\exp\left\{-n\left(\frac{x-a }{\omega}\right)^{-\bar{\beta}}\right\} 
\\ &\quad \text{for } \bar{\beta}>0 
\\ \notag
& = \text{GenFisherTippett}(x|a, \omega, n, -\bar{\beta})
\\ \notag
&= \text{Amoroso}(x| a , {\omega}/{n^{\frac{1}{\beta} }},n,-\bar{\beta}),
\end{align}
The limiting distribution of the $n$th largest value of a large number identically distributed random variables whose moments are not all finite and are bounded from below by $a$.  (If the shape parameter $\omega$ is negative then minimum rather than maxima.)

\dist{Fr\'{e}chet} (extreme value type II , Fisher-Tippett type II, Gumbel type II, inverse Weibull) distribution~\cite{Frechet1927,Gumbel1958}:
\begin{align}
\label{Frechet}
\text{Fr\'{e}chet}(x | a , \omega, \bar{\beta}) 
&=	\frac{\bar{\beta}}{| \omega |} \left(\frac{x-a }{\omega}\right)^{-\bar{\beta}-1} 
\exp\left\{-\left(\frac{x-a }{\omega}\right)^{-\bar{\beta}}\right\}
\\ \notag &\quad \text{for } \bar{\beta}>0 
\\  \notag
& = \text{FisherTippett}(x| a, \omega, -\bar{\beta}) 
\\ \notag 
&= \text{Amoroso}(x| a , \omega,1,-\bar{\beta})\notag
\end{align}
The limiting distribution of the largest of a large number identically distributed random variables whose moments are not all finite and are bounded from below by $a$.  (If the shape parameter $\omega$ is negative then minimum rather than maxima.)
Special cases of the Fr\'{e}chet  distribution include the inverse exponential ($\bar{\beta}=1$) and inverse Rayleigh ($\bar{\beta}=2$) distributions.

\subsection{Properties}
\begin{align*}
\text{support} \quad & x \geq a &  \theta > 0
\\
&   x\leq a  &  \theta < 0 
 \\
\text{cdf} \quad 
&
 1-Q(\alpha, \left(\tfrac{x - a }{\theta}\right)^{\beta}) 
& \tfrac{\beta}{\theta}>0
\\
& 
 Q(\alpha, \left(\tfrac{x - a }{\theta}\right)^{\beta}) 
& \tfrac{\beta}{\theta}<0
\\
\text{mode} \quad&   a+ \theta (\alpha-\tfrac{1}{\beta})^{\frac{1}{\beta}} 
& \alpha \beta  \geq 1
\\ & a & \alpha \beta  \le 1
\\
\text{std. moments}\quad &  \frac{\Gamma(\alpha+\frac{r}{\beta})}{\Gamma(\alpha)}  \quad ( a=0, \theta=1) 
& \alpha + \tfrac{n}{\beta} \geq 0
\\
\text{mean} \quad& a  + \theta \frac{\Gamma(\alpha+\frac{1}{\beta})}{\Gamma(\alpha)} 
& \alpha + \tfrac{1}{\beta} \geq 0
\\
\text{variance}  \quad&   \theta^2 \left[  \frac{\Gamma(\alpha+\frac{2}{\beta})}{\Gamma(\alpha)}  - 
\frac{\Gamma(\alpha+\frac{1}{\beta})^2}{\Gamma(\alpha)^2}    \right]
& \alpha + \tfrac{2}{\beta} \geq 0
 \\
\text{entropy}  \quad& 
\ln \frac{\theta \Gamma(\alpha)}{|\beta|} +\alpha + \left( \frac{1}{\beta} - \alpha\right) \psi(\alpha) &
\text{\cite{Dadpay2007}}
\end{align*}
Here, cdf is the cumulative distribution function, 
$Q(\alpha,x)= \Gamma(\alpha,x) /\Gamma(\alpha)$ is the regularized gamma function~\cite{Abramowitz1965},
$\Gamma(\alpha,x) = \int_x^{\infty} t^{\alpha-1} e^{-t} dt$ is the incomplete gamma integral~\cite{Abramowitz1965},
and
$\psi(x)=\frac{d}{dx} \ln \Gamma(x)$ is the digamma function~\cite{Abramowitz1965}, the logarithmic derivative of the gamma function.
 Important special cases and limits include
$\Gamma(\tfrac{1}{2})=\sqrt{\pi}$, $\Gamma(\tfrac{1}{2}, x) = \sqrt{\pi}\ \text{erfc}(\sqrt{x})$ and $\Gamma(1,x) = \exp(-x)$.
The derivative of the regularized gamma function is $\tfrac{d}{dx} Q(\alpha,x) = -\tfrac{1}{\Gamma(\alpha)} x^{\alpha-1} e^{-x}$.

The profile of the Amoroso distribution is bell shaped for $\alpha\beta\geq1$, and otherwise L- or J- shaped with the  mode at the boundary. The moments are undefined if the side conditions are not satisfied. Expressions for skew and kurtosis are not simple, but can be deduced from the moments if necessary.

The Amoroso distribution can be obtained from the standard gamma distribution (\ref{StdGamma}) with the change of variables,  $x \mapsto \left(\tfrac{x-a}{\theta}\right)^\beta$. Therefore, Amoroso random numbers can be obtained by sampling from the standard gamma distribution, for instance using the Marsaglia-Tsang fast gamma method~\cite{Marsaglia2001} and applying the appropriate  transformation~\cite{Knuth1997}.

\pagebreak[4]
\section{Log-gamma distributions}
\label{SecLogGamma}

The {\bf log-gamma} (Coale-McNeil, gamma-exponential) distribution~\cite{Bartlett1946,Prentice1974,Johnson1995} is  a three parameter, continuous, univariate, unimodal probability density with infinite range. The functional form in the most straightforward parameterization is
\begin{align}
\label{LogGamma}
\text{LogGamma}(x|\nu, \lambda, \alpha) 
&=
\frac{1}{ \Gamma(\alpha) |\lambda|}  \exp\left\{ \alpha \left(\frac{x-\nu}{\lambda}\right) - \exp\left(\frac{x-\nu}{\lambda}\right)  \right\} 
\\ \notag
& \qquad \text{for } x,\ \nu,\ \lambda,\ \alpha,\   \text{in } \mathbb{R}, 
\ \alpha>0, \ 
\\ \notag
& \qquad \text{support } -\infty \leq x \leq \infty
\notag
\end{align}
The three real parameters consist of a location parameter $\nu$, a scale parameter~$\lambda$, and a shape parameter $\alpha$, which is inherited directly from the Amoroso distribution.

The name ``log-gamma'' arises because the standard log-gamma distribution is the logarithmic transform of the standard gamma distribution
\begin{align*}
\text{StdLogGamma}(\alpha)  &\sim \ln\Big( \text{StdGamma}(\alpha) \Big)
\\
\text{LogGamma}(\nu, \lambda, \alpha)  &\sim \ln\Big( \text{Amoroso}(0, e^{\nu},\alpha, \tfrac{1}{\lambda}) \Big)
\end{align*}
Note that this naming convention is the opposite of that used for the log-normal distribution (\ref{LogNormal}). The name  ``log-gamma''  has also been used for the antilog transform of the generalized gamma distribution, which leads to the unit-gamma distribution~\cite{Gupta2004}.

The log-gamma distribution is a limit of the Amoroso distribution (\ref{Amoroso}),  and itself  has a number of important limits and special cases (Table II), which we will discuss below. 
\begin{align}
&\text{LogGamma}(x|\nu, \lambda, \alpha) 
\\ \notag & = \lim_{\beta\rightarrow\infty}   
\frac{1}{\Gamma(\alpha) |\lambda|} 
\left(1+ \frac{1}{\beta}\frac{x-\nu }{\lambda}\right)^{\alpha\beta-1}
\!\!\exp\left\{  - \left(1+ \frac{1}{\beta}\frac{x-\nu }{\lambda}\right)^{\beta} \right\}
\\ \notag &= \lim_{\beta\rightarrow\infty} \text{Amoroso}(\nu-\beta \lambda,\beta \lambda, \alpha,\beta)
\notag
\end{align}
Recall that $\lim_{\beta\rightarrow\infty} (1+ \tfrac{x}{\beta})^{\beta} = \exp(x) $.

\begin{table}
\caption{The log-gamma family of distributions.}
\begin{align*}
&\text{LogGamma}(x|\nu, \lambda, \alpha) 
\\&=
\frac{1}{ \Gamma(\alpha) | \lambda |}  \exp\left\{ \alpha \left(\frac{x-\nu}{\lambda}\right) - \exp\left(\frac{x-\nu}{\lambda}\right)  \right\} \notag
\\ \notag
& \qquad \text{for } x,\ \nu,\ \lambda,\ \alpha,\   \text{in } \mathbb{R}, 
\ \alpha>0, \ 
\\ \notag
& \qquad \text{support } -\infty \leq x \leq \infty
\notag
\end{align*}

\begin{tabular}{llcccl}
(\ref{LogGamma}) & log-gamma &  $\nu$ & $\lambda$ & $\alpha$ 
\\ \hline
(\ref{StdLogGamma}) & standard log-gamma &  $0$ & $1$ & $\alpha$  \\
(\ref{LogChiSqr}) & log-chi-square &$\ln 2 $ & $1$ & $\tfrac{k}{2}$ \\ 
(\ref{GenGumbel}) &generalized Gumbel    & . & . & $n$ &  \\
(\ref{Gumbel}) &Gumbel  &  . & . & 1 &  \\
(\ref{BHP}) &BHP   & . & . & $\frac{\pi}{2}$ &  \\
(\ref{StdGumbel}) &standard Gumbel &   0 & -1 & 1 & \\ 
\\
 & \underline{Limits} \\
(\ref{Normal}) & normal  & . &  . & $\alpha$& $\lim_{\alpha\rightarrow\infty}$  \\

\end{tabular}
\end{table}

\begin{figure}[t]
\begin{center}
\includegraphics[scale=1]{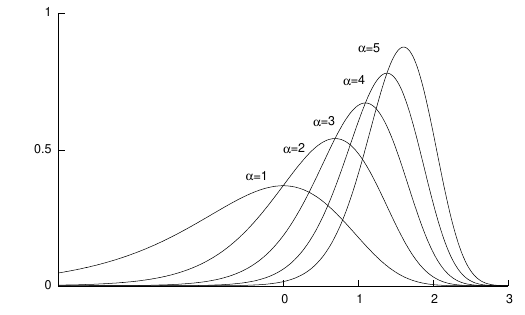}
\end{center}
\caption{$\text{\rm LogGamma}(x|0,1,\alpha)$}
\end{figure}

\subsection{Special cases}

\dist{Standard log-gamma} distribution:
\begin{align}
\label{StdLogGamma}
\text{StdLogGamma}(x| \alpha) 
=&
\frac{1}{\Gamma(\alpha) }  \exp\left\{ \alpha x - \exp (x)  \right\} 
\\ \notag =& \text{LogGamma}(x|0,1,\alpha)
\end{align} 
The log-gamma distribution with zero location and unit scale.

\dist{Log-chi-square} distribution~\cite{Lee2012}:
\begin{align}
\label{LogChiSqr}
\text{LogChiSqr}(x|k) 
 &= \notag
\frac{1}{2^{\frac{k}{2}} \Gamma(\frac{k}{2})}  \exp\left\{ \frac{k}{2} x - \frac{1}{2} \exp(x)  \right\} 
\\ & \qquad \text{ for positive integer } k
\\
&= \text{LogGamma}(x|\ln 2, 1 ,\tfrac{k}{2}) 
\notag
\end{align}
The logarithmic transform of the chi-square distribution~(\ref{ChiSqr}).

\dist{Generalized Gumbel} distribution~\cite{Gumbel1958,Johnson1995}: 
\begin{align}
\label{GenGumbel}
\text{GenGumbel}(x|u,\bar{\lambda},n) 
 &= \notag
\frac{n^n}{\Gamma(n) |\bar{\lambda}|}  \exp\left\{ - n \left(\frac{x-u}{\bar{\lambda}}\right) - n \exp\left(- \frac{x-u}{\bar{\lambda}}\right)  \right\} 
\\ & \qquad \text{ for positive integer } n
\\
&= \text{LogGamma}(x|u+\bar{\lambda} \ln n,-\bar{\lambda},n) 
\notag
\end{align}
The limiting distribution of the $n$th largest value of a large number of unbounded identically distributed random variables whose probability distribution has an exponentially decaying tail.

\dist{Gumbel}  (Fisher-Tippett type I, Fisher-Tippett-Gumbel, FTG, Gumbel-Fisher-Tippett, log-Weibull, extreme value (type I),  doubly exponential, double exponential) distribution~\cite{Fisher1928,Gumbel1958, Johnson1995}:
\begin{align}
\label{Gumbel}
\text{Gumbel}(x|u,\bar{\lambda}) 
&=
\frac{1}{|\bar{\lambda}|}  \exp\left\{ -\left(\frac{x-u}{\bar{\lambda}}\right) - \exp\left(-\frac{x-u}{\bar{\lambda}}\right)  \right\} 
\\
&= \text{LogGamma}(x|u,-\bar{\lambda},1) 
\notag
\end{align}
This is the asymptotic extreme value distribution for variables of ``exponential type'', unbounded with finite moments~\cite{Gumbel1958}.
With positive scale $\bar{\lambda}>0$, this is an extreme value distribution of the maximum, with negative scale $\bar{\lambda}<0$ ($\lambda>0$) an extreme value distribution of the minimum. Note that the Gumbel is sometimes defined with the negative of the scale used here.

Note that the term ``double exponential distribution'' can refer to either the Gumbel or Laplace~\cite{Johnson1995} distributions.
The Gompertz distribution is a truncated  Gumbel distribution~\cite{Johnson1995}. 

\dist{Standard Gumbel} (Gumbel) distribution~\cite{Gumbel1958}:
\begin{align}
\label{StdGumbel}
\text{StdGumbel}(x) 
=&
  \exp\left\{- x - e^{-x} \right\}  \\
=& \text{LogGamma}(x|0,-1,1)  \notag
\end{align}
The Gumbel distribution with zero location and a unit scale.

\dist{BHP} (Bramwell-Holdsworth-Pinton) distribution~\cite{Bramwell1998}:
\begin{align}
\label{BHP}
\text{BHP}(x|\nu,\lambda) 
&=
\frac{1}{\Gamma(\tfrac{\pi}{2}) |\lambda|}  \exp\left\{ \frac{\pi}{2}\left(\frac{x-\nu}{\lambda}\right) -  \exp\left(\frac{x-\nu}{\lambda}\right)  \right\} 
\\
&= \text{LogGamma}(x|\nu,\lambda,\frac{\pi}{2})  
\notag
\end{align}
Proposed as a model of rare fluctuations in turbulence and other correlated systems.

\subsection{Properties}

\begin{align*}
\text{support}\quad & -\infty \leq x \leq +\infty & \\
\text{cdf} \quad &  1- Q(\alpha, e^{\frac{x-\nu}{\lambda} }) & \text{ for } \lambda>0 
\\ & Q(\alpha, e^{\frac{x-\nu}{\lambda} }) & \text{ for } \lambda<0 
\\
\text{cgf} \quad & \nu t + \ln \frac{\Gamma(\alpha+\lambda t)}{\Gamma(\alpha)}  & \text{\cite{Johnson1995}}
\\
\text{mode} \quad &  \nu -\lambda \ln \alpha
\\
\text{mean} \quad & \nu+ \lambda \psi(\alpha) & 
\\
\text{variance}  \quad &   \lambda^2 \psi_1(\alpha) & 
\\
\text{skew}  \quad &  \text{sgn}(\lambda) \psi_2(\alpha) /  \psi_1(\alpha)^{3/2} 
\\
\text{kurtosis}  \quad &     \psi_3(\alpha) /  \psi_1(\alpha)^{2}
 \\
\text{entropy}  \quad & \ln \Gamma(\alpha) |\lambda|  - \alpha \psi(\alpha) + \alpha & 
\end{align*}
Here, cdf is the cumulative distribution function, cgf the cumulant generating function, $\ln \text{E}[ \exp( t X)]$, $\text{sgn}(z)$ is the sign function, $\psi_n(z)= \tfrac{d^n}{dz^n} \ln \Gamma(z)$ is the polygamma function and  $\psi(z)\equiv\psi_0(z)$ is the digamma function.

\section{Miscellaneous limits}

\label{SecMisc}
 
\dist{Log-normal} (Galton, Galton-McAlister, antilog-normal,  logarithmic-normal, logarithmico-normal, Cobb-Douglas, $\Lambda$) distribution~\cite{Galton1879, McAlister1879, Johnson1994}: 
 \begin{align}
 \label{LogNormal}
 \text{LogNormal}(x|a, \vartheta,\sigma) 
&= \frac{\vartheta}{\sqrt{2\pi \sigma^2}} \left(\frac{x-a}{\vartheta}\right)^{-1} \exp\left\{-\frac{1}{2\sigma^2} \left( \ln \frac{x-a}{\vartheta} \right)^2 \right\}
\\ \notag & 
=  \lim_{\beta\rightarrow0} \text{Amoroso}(x| a, \vartheta (\beta\sigma)^{2/\beta} ,1/(\beta\sigma)^2, \beta)
\end{align}
The log-normal distribution is a limiting form of the gamma family.  To see this, make the requisite substitutions, 
\begin{align*}
& \text{Amoroso}(x| a, \vartheta (\beta\sigma)^{2/\beta} ,1/(\beta\sigma)^2, \beta)  \\
&\propto
 \left(\frac{x-a}{\vartheta}\right)^{-1} 
\exp\left\{
\frac{1}{\beta\sigma^2} \ln \frac{x-a}{\vartheta}  
- \frac{1}{\beta^2\sigma^2} \exp \left[ \beta \ln \frac{x-a}{\vartheta} \right]
\right\}
\notag
\end{align*}
and in the limit ${\beta\rightarrow0}$  expand the second exponential to second order in $\beta$.

With $a=0$, $\vartheta=1$, $\sigma=1$  we obtain the  {\bf standard log-normal} (Gibrat) distribution~\cite{Gibrat1931}.
The {\bf two-parameter lognormal} distribution ($a=0$) arises from the multiplicative version of the central limits theorem: When the sum of independent random variables limits to normal, the product of those random variables limits to log-normal.
The log-normal distribution maps to the normal distribution with the transformation $x\mapsto \exp(x)$.
\begin{align*}
\text{LogNormal}(a, \vartheta, \sigma) &\sim \exp\Big(\text{Normal}(\ln \vartheta,\sigma)\Big) +a
\end{align*}

 \dist{Normal} (Gauss, Gaussian, bell curve, Laplace-Gauss, de Moivre, error, Laplace's second law of error, law of error) distribution~\cite{Moivre1733,Johnson1994}:
\begin{align}
\label{Normal}
\text{Normal}(x|\mu,\sigma) 
&=
\frac{1}{\sqrt{2 \pi \sigma^2}}  \exp\left\{ - \frac{( x-\mu)^2}{2\sigma^2} \right\}
 \\ \notag
&= \lim_{\alpha\rightarrow\infty} \text{Amoroso}(x| \mu- \sigma\sqrt{\alpha}, \sigma/\sqrt{\alpha}, \alpha, 1)
 \\ \notag
&= \lim_{\alpha\rightarrow\infty} \text{LogGamma}(x| \mu-\sigma\sqrt{\alpha}\ln\alpha, \sigma\sqrt{\alpha}, \alpha) 
\end{align}
With $\mu=0$ and $\sigma = 1/ \sqrt{2} h$ we obtain the {\bf error function} distribution, and
with $\mu=0$ and $\sigma=1$ we obtain the {\bf standard normal} ($\Phi$, $z$, unit normal)  distribution. In the limit that $\sigma\rightarrow\infty$ we obtain an unbounded {\bf uniform} (flat distribution, and in the limit $\sigma\rightarrow0$ we obtain a {\bf delta} (degenerate) distribution. 

The normal distribution is a limit of the Amoroso~\cite{Johnson1994} and log-gamma distributions~\cite{Prentice1974}. For Amoroso, make the requisite substitutions,
\begin{align*}
& \text{Amoroso}(x|  \mu- \sigma\sqrt{\alpha}, \sigma/\sqrt{\alpha}, \alpha,1)  \\
&\propto \
\exp\left\{ -\sqrt{\alpha} \frac{x-\mu}{\sigma} + (\alpha-1) \ln \left( 1+ \frac{1}{\sqrt{\alpha}} \frac{x-\mu}{\sigma}  \right) \right\}
\end{align*}
and expand the logarithm as $\ln(1+x) = x - \tfrac{x^2}{2} + O(x^3)$.

\dist{Power-law} (Pearson type XI, fractal) distribution~\cite{Pearson1916}: 
\begin{align}
\label{PowerLaw}
\text{PowerLaw}(x|p) 
&\propto   \frac{1}{(x-a)^p}
\\&=  \lim_{\alpha\rightarrow\infty} \text{Amoroso}(a,\theta,\alpha,(1-p)/\alpha) \notag
\end{align}
Improper (unnormalizable) power law distributions are obtained as a limit of the gamma distribution family.
 If $p=0$ we obtain the {\bf half-uniform} distribution over the positive numbers; if $p=1$ we obtain {\bf Jeffreys} distribution~\cite{Jeffreys1948}, used as an uninformative  prior in Bayesian probability~\cite{Jaynes2003}. 

~

~

~

\paragraph*{Acknowledgments:}
I am grateful to  David Sivak, Edward E. Ayoub and Francis J. O'Brien for spotting various typos  and errors,  and, as always, to Avery Brooks for many insightful observations. In curating this collection of distributions, I have benefited mightily from Johnson, Kotz, and Balakrishnan's monumental compendiums~\cite{Johnson1994, Johnson1995}, Eric Weisstein's MathWorld, and the myriad pseudo-anonymous contributors to Wikipedia. An extended enumeration of simple probability distribution families can be found at \href{http://threeplusone.com/gud}{threeplusone.com/gud}.

\appendix

\section{Index of distributions}


\paragraph*{generalized-X}  The only consistent meaning is that distribution ``X'' is a special case of the distribution ``generalized-X''.  In practice, often means ``add a shape parameter''. 

\paragraph*{standard-X} The distribution ``X'' with the location parameter set to 0, scale to 1, and often the Weibull shape parameter $\beta$ to 1. Not to be confused with {\it standardized} which generally indicates zero mean and unit variance. 

\paragraph*{shifted-X}  (or translated) A distribution with an additional location parameter. 

\paragraph*{scaled-X} (or scale-X) A distribution with an additional scale parameter. 

\paragraph*{inverse-X} (Occasionally inverted-X, reciprocal-X, or negative-X) Generally labels the transformed distribution with $x\mapsto\tfrac{1}{x}$, or more generally the distribution with the Weibull shape parameter negated, $\beta\rightarrow-\beta$. An exception is the inverse normal distribution~\cite{Johnson1994}. 

\paragraph*{log-X}  If $x$ follows distribution X then either $y=\ln x$ (e.g. log-normal) or $y = e^x$ (e.g. log-gamma) follows log-X.  This ambiguity arrises because although the second convention is more logical, the log-normal convention has historical precedence.

\paragraph*{reversed-X} (Occasionally negative-X)  The scale is negated. 
\paragraph*{X of the Nth kind} See ``X type N''.
~

~



\noindent
Distribution\hfill Synonym or Equation \\
\\
%
$\chi$ 						\dotfill  chi\\
$\chi^2$ 						\dotfill  chi-square \\
$\Gamma$ 					\dotfill  gamma \\
$\Lambda$					\dotfill log-normal \\ 
$\Phi$						\dotfill standard normal \\
antilog-normal 					\dotfill log-normal \\
Amoroso						\dotfill(\ref{Amoroso})\\
bell curve 						\dotfill  normal \\
BHP							\dotfill (\ref{BHP})\\
Bramwell-Holdsworth-Pinton 		\dotfill  BHP \\
%
chi							\dotfill (\ref{Chi})\\
chi-square					\dotfill (\ref{ChiSqr})\\
Coale-McNeil 					\dotfill  generalized log-gamma \\  
Cobb-Douglas					\dotfill	log-normal \\
%
de Moivre						\dotfill normal \\
degenerate					\dotfill delta \\
delta  						\dotfill  See normal (\ref{Normal})  \\ 
doubly exponential 				\dotfill  Gumbel \\
double exponential 				\dotfill  Gumbel or Laplace \\
Erlang						\dotfill  See gamma (\ref{Gamma}) \\
error 						\dotfill  normal \\
error function					\dotfill See normal (\ref{Normal}) \\
exponential 					\dotfill (\ref{Exp}) \\ 
extreme value					\dotfill Gumbel \\
extreme value type N 			\dotfill  Fisher-Tippett type N \\
%
Fisher-Tippett					\dotfill  (\ref{FisherTippett}) \\
Fisher-Tippett type I				\dotfill  Gumbel \\
Fisher-Tippett type II				\dotfill   Fr\'{e}chet  \\
Fisher-Tippett type III			\dotfill  Weibull \\
Fisher-Tippett-Gumbel			\dotfill  Gumbel  \\
fractal 						\dotfill  power law \\
flat 							\dotfill  uniform \\
Fr\'{e}chet 					\dotfill (\ref{Frechet}) \\
FTG							\dotfill  Fisher-Tippett-Gumbel \\
Galton						\dotfill log-normal \\
Galton-McAlister				\dotfill log-normal \\
gamma						\dotfill (\ref{Gamma}) \\
gamma-exponential				\dotfill log-gamma \\
Gaussian 						\dotfill  normal \\
Gauss 						\dotfill  normal \\
generalized gamma 				\dotfill Stacy or Amoroso \\
generalized inverse gamma   		\dotfill See Stacy (\ref{Stacy}) \\
generalized Gumbel 			\dotfill (\ref{GenGumbel}) \\
generalized extreme value 		\dotfill Fisher-Tippett \\
generalized Fisher-Tippett		\dotfill (\ref{GenFisherTippett}) \\
generalized Fr\'{e}chet 			\dotfill  (\ref{GenFrechet}) \\
generalized inverse gamma 		\dotfill   generalized gamma \\
generalized normal	 			\dotfill  Nakagami \\
generalized Rayleigh	 		\dotfill   scaled chi\\
generalized semi-normal			\dotfill Stacy \\   
generalized Weibull				\dotfill (\ref{GenWeibull}) \\
GEV 							\dotfill  generalized extreme value \\
Gibrat						\dotfill standard log-normal \\ 
Gumbel 						\dotfill (\ref{Gumbel}) \\
Gumbel-Fisher-Tippett			\dotfill  Gumbel \\ 
Gumbel type N 					\dotfill  Fisher-Tippett type N \\
%
half-normal 					\dotfill (\ref{HalfNormal})\\
half-uniform 					\dotfill See power law (\ref{PowerLaw}) \\
hydrograph					\dotfill Stacy \\
hyper gamma   					\dotfill  Stacy\\
%
inverse chi  					\dotfill (\ref{InvChi})\\
inverse chi-square 				\dotfill (\ref{InvChiSqr})\\
inverse exponential 				\dotfill  (\ref{InvExp})\\
inverse gamma 				\dotfill (\ref{InvGamma})\\
inverse Rayleigh  				\dotfill (\ref{InvRayleigh})\\
inverse Weibull 				\dotfill  Fr\'{e}chet  \\
Jeffreys  						\dotfill See power law (\ref{PowerLaw})\\
%
%
Laplace's second law of error  		\dotfill normal \\
Laplace-Gauss					\dotfill normal \\
law of error					\dotfill normal \\
Leonard hydrograph			\dotfill Stacy \\
L\'{e}vy 						\dotfill (\ref{Levy})\\
log-chi-square					\dotfill (\ref{LogChiSqr}) \\
log-gamma 					\dotfill 	(\ref{LogGamma}) \\
log-normal					\dotfill 	(\ref{LogNormal})	\\
log-normal, two parameter		\dotfill	See log-normal (\ref{LogNormal}) \\
log-Weibull 					\dotfill	Gumbel \\
logarithmic-normal				\dotfill log-normal \\
logarithmico-normal				\dotfill log-normal \\
%
Maxwell 						\dotfill (\ref{Maxwell}) \\
Maxwell-Boltzmann 				\dotfill  Maxwell \\ 
Maxwell speed					\dotfill Maxwell \\
m-Erlang						\dotfill Erlang \\
Nakagami 					\dotfill (\ref{Nakagami}) \\
Nakagami-m 					\dotfill  Nakagami \\    
negative exponential			\dotfill exponential \\
normal 						\dotfill (\ref{Normal})\\
Nukiyama-Tanasawa 			\dotfill  generalized gamma \\	
one-sided normal  				\dotfill  half-normal\\
%
Pearson type III   				\dotfill (\ref{PearsonIII}) \\
Pearson type V    				\dotfill (\ref{PearsonV}) \\
Pearson type X    				\dotfill  exponential\\
Pearson type XI    				\dotfill  power law\\				
positive definite normal 			\dotfill  half-normal \\ 
power law 					\dotfill (\ref{PowerLaw}) \\
pseudo-Weibull				\dotfill (\ref{PseudoWeibull})\\
%
%
Rayleigh  						\dotfill (\ref{Rayleigh}) \\
Rosin-Rammler 				\dotfill  Weibull\\				
Rosin-Rammler-Weibull 			\dotfill  Weibull\\
%
scaled chi						\dotfill (\ref{ScaledChi}) \\
scaled chi-square 				\dotfill  (\ref{ScaledChiSqr}) \\
scaled inverse chi 				\dotfill (\ref{ScaledInvChi})\\
scaled inverse chi-square 		\dotfill  (\ref{ScaledInvChiSqr})\\ 
semi-normal 					\dotfill  half-normal \\ 
shifted exponential 				\dotfill  (\ref{ShiftExp}) \\
Stacy 						\dotfill   (\ref{Stacy}) \\
Stacy-Mihram					\dotfill Amoroso\\
standard Amoroso				\dotfill standard gamma \\
standard exponential			\dotfill See exponential (\ref{Exp}) \\
standard gamma 				\dotfill (\ref{StdGamma}) \\
standard Gumbel 				\dotfill (\ref{StdGumbel}) \\
standard log-gamma			\dotfill (\ref{StdLogGamma}) \\
standard log-normal				\dotfill See log-normal (\ref{LogNormal}) \\
standard normal 				\dotfill See normal (\ref{Normal}) \\
stretched exponential 			\dotfill   (\ref{StretchedExp}) \\
%
transformed gamma				\dotfill Stacy \\
uniform  						\dotfill  See normal (\ref{Normal})  \\ 
unit normal					\dotfill standard normal \\
van der Waals profile 			\dotfill  L\'{e}vy \\
Vienna 						\dotfill  Wien \\
Vinci 						\dotfill  inverse gamma \\
von Mises extreme value			\dotfill Fisher-Tippett \\
von Mises-Jenkinson 			\dotfill  Fisher-Tippett \\
waiting time					\dotfill exponential \\
Weibull 						\dotfill (\ref{Weibull}) \\
Weibull-Gnedenko				\dotfill Weibull \\
Wien							\dotfill  See gamma (\ref{Gamma})\\
Wilson-Hilferty					\dotfill (\ref{WilsonHilferty}) \\
$z$							\dotfill standard normal \\




\end{document}